\documentclass[12pt]{amsart}
\usepackage{amssymb} 

\newtheorem{thm}{Theorem}[section]

\theoremstyle{definition}

\newtheorem{rem}[thm]{Remark}
\numberwithin{equation}{section}

\begin{document}
\title[]{Finite $p$-groups of class 2 have noninner
automorphisms of order $p$}%
\author[A. Abdollahi]{A. Abdollahi}%
\address{Department of Mathematics, University of Isfahan, Isfahan 81746-73441, Iran; and Institute for Studies in Theoretical
Physics and Mathematics (IPM). }%
\email{a.abdollahi@math.ui.ac.ir}%
\thanks{This research was in part supported by a grant from IPM (No. 85200032). }
\subjclass{20D45;20E36}
\keywords{Automorphisms of $p$-groups; nilpotent groups of class 2; noninner automorphisms}%
\dedicatory{Dedicated to Professor Aliakbar Mohammadi Hassanabadi }%
\begin{abstract} We prove
that for any prime number $p$, every finite non-abelian $p$-group
 $G$ of class 2 has a noninner automorphism of order $p$ leaving either the
Frattini subgroup $\Phi(G)$ or  $\Omega_1(Z(G))$ elementwise
fixed.
\end{abstract}
\maketitle
\section{\bf Introduction}
Let $p$ be a prime number and $G$ be a non-abelian finite
$p$-group. A longstanding conjecture asserts that  $G$ admits a
noninner automorphism of order $p$ (see also Problem 4.13 of
\cite{Kbook}). By a famous result of W. Gasch\"utz \cite{G},
noninner automorphisms of $G$ of $p$-power order exist. M.
Deaconescu and G. Silberberg \cite{DS} reduced the verification
of the conjecture to the case in which $C_G(Z(\Phi(G)))=\Phi(G)$.
H. Liebeck \cite{L} has shown that finite $p$-groups of class 2
with $p>2$ must have a noninner automorphism of order $p$ fixing
the Frattini subgroup elementwise.  It follows from a
cohomological result of P. Schmid \cite{S} that the conjecture is
true whenever $G$ is regular. Here we show the validity of  the
conjecture when $G$ is nilpotent of class 2. In fact we prove that \\

\noindent {\bf Theorem.} For any prime number $p$, every finite
non-abelian $p$-group  $G$ of class 2 has a noninner automorphism
of order $p$ leaving either the Frattini subgroup $\Phi(G)$ or
$\Omega_1(Z(G))$ elementwise fixed.\\

The unexplained notation is standard and follows that of
Gorenstein \cite{Gor}.

\section{\bf Preliminaries}
We use the following  facts in the proof of the Theorem.
\begin{rem}\label{G'}
If $G$ is a  group whose derived subgroup $G'$ is a finite cyclic
$p$-group for some prime $p$, then $G'=\langle [a,b]\rangle$ for
some $a,b\in G$. Since $G'$ is generated by commutators $[x,y]$
($x,y\in G$) whose orders are $p$-powers and $G'$ is abelian,
$exp(G')=\max\{ |[x,y]| \;:\; x,y\in G\}$. But $G'$ is a finite
cyclic group and so $exp(G')=|G'|$. Hence $G'$ is generated by
one of the elements of the set $\{ [x,y] \;:\; x,y\in G\}$.
\end{rem}
\begin{rem}\label{H}
Let $G$ be a finite nilpotent group of class 2 such that
$G'=\langle [a,b] \rangle$ for some $a,b \in G$. Then by a
well-know argument (e.g., see the proof of Lemma 1 of \cite{AY})
we have $G=\left<a,b\right>C_G(\langle a,b\rangle)$. We give it
here for the reader's convenience: for any $x\in G$, we have
$[a,x]=[a,b]^s$ and $[b,x]=[a,b]^t$ for some integers $s,t$. Then
$[a,b^{-s}a^tx]=1$ and $[b,b^{-s}a^tx]=1$. Hence $b^{-s}a^tx\in
C_G(\langle a,b\rangle)$ and so $G=\left<a,b\right>C_G(\langle
a,b\rangle)$.
\end{rem}
\begin{rem}\label{iden}
Let $G$ be a nilpotent group of class 2, $x,y\in G$ and $k>0$ be
an integer. Then since $[y,x]=y^{-1}x^{-1}yx\in Z(G)$, it is easy
to see by induction on $k$ that
$(xy)^k=x^ky^k[y,x]^{\frac{k(k-1)}{2}}$. Also we have
$[x,y]^m=[x^m,y]=[x,y^m]$ for all integers $m$.
\end{rem}
 We shall make frequent use of Remark
\ref{iden} without reference in the proof of the Theorem.
Especially we use it in such a sample situation: if we know that
$x$ and $y$ are two elements in a nilpotent $2$-group of class 2,
 $m\in\mathbb{Z}$ and $n\in\mathbb{N}$ such that $|[x,y]|=2^n$ and
$x^{m2^n}=y^{-2^n}$, then by Remark \ref{iden} and the
hypothesis  we have
$$(x^my)^{2^n}=x^{m2^n}y^{2^n}[y,x^m]^{2^{n-1}(2^n-1)}=[y,x]^{m2^{n-1}(2^n-1)}.$$
Since $[x,y]=[x,x^my]$ and $|[x,y]|=2^n$, we have that
$(x^my)^{2^{n-1}}\not=1$ and so $2^n \mid |x^my|$. It follows
that $|x^my|=2^{n+1}$, if $m$ is odd, and $|x^my|=2^n$, if $m$ is
even.
\begin{rem}\label{Lie}
Let $G$ be a finite $p$-group of class 2. If $G$ has  no noninner
automorphism of order $p$ leaving $\Phi(G)$ elementwise fixed,
then $Z(G)$ must be cyclic. In fact  by the part (a) of the proof
of \cite[Theorem 1]{L}, we have $G'$  is cyclic. Now if  $Z(G)$ is
not cyclic, then $\Omega_1(Z(G))$ is not cyclic and so
$\Omega_1(Z(G))\nleq G'$. Now take an element  $z\in
\Omega_1(Z(G))\backslash G'$, a maximal subgroup $M$ of $G$ and
$g\in G\backslash M$. Then the map $\alpha$ on $G$ defined by
$(mg^i)^\alpha=mg^iz^i$ for all $m\in M$ and integers $i$, is a
noninner automorphism of order $p$ leaving $M$ (and so $\Phi(G)$)
elementwise fixed, a contradiction.\\
Note that if $Z(G)=\Phi(G)$, then one may replace the latter
argument by the part (iv) of Lemma 2 of \cite{L}.
\end{rem}
\begin{rem}\label{cent}
Let $G$ be a group and $H,K$ be subgroups of $G$ such that $G=HK$
and $[H,K]=1$. If there exists a noninner automorphism $\varphi$
of order $p$ in $Aut(H)$ leaving $Z(H)$ elementwise fixed, then
the map $\beta$ on $G$ defined by $(hk)^\beta=h^\varphi k$ for
all $h\in H$ and $k\in K$ is a noninner automorphism of $G$ of
order $p$ leaving $Z(G)$ elementwise fixed. It is enough to show
that $\beta$ is well-defined and this  can be easily seen, because
$x^\varphi=x$ for all $x\in H \cap K=Z(H)$, by hypothesis.
\end{rem}
\section{\bf Proof of the Theorem}
By the main results of \cite{DS} and \cite{L}, we may assume that
$\Phi(G)=C_G(Z(\Phi(G)))$ and $p=2$. By Remark \ref{Lie}, we may
further  assume that $Z(G)$ is cyclic.
Now Remark \ref{G'} implies that there exist elements $a,b\in G$
such that $G'=\left<[a,b]\right>$.  Let $H=\left<a,b\right>$.
Then it follows from  Remark \ref{H}  that $G=HC_G(H)$ and by
Remark \ref{cent} it is enough to construct
 a noninner automorphism $\varphi$ of $H$ of order $2$ leaving $Z(H)$ elementwise fixed. \\
Note that $|G'|=|H'|=|[a,b]|=2^n$ for some integer $n>0$. Since
$G'$ is cyclic and $G'\leq Z(G)$,
$$\text{exp}(\frac{G}{Z(G)})=\text{exp}(\frac{H}{Z(H)})=2^n,$$ which implies   that
$Z(H)=\left<a^{2^n},b^{2^n},[a,b]\right>\leq Z(G)$. If $n=1$,
then $\Phi(G)=G^2\leq Z(G)$. Since $\Phi(G)=C_G(Z(\Phi(G)))$, we
have $G=\Phi(G)$, which is impossible. Therefore $n\geq 2$. Since
$Z(H)$ is cyclic, either $a^{2^ni}=b^{2^n}$ or $a^{2^n}=b^{2^ni}$
for some integer $i$.  Suppose that $a^{2^ni}=b^{2^n}$. If $i$ is
even, then $|a^{-i}b|=2^n$ and $(a^{-i}b)^{2^{n-1}}\notin Z(H)$,
as $[a,b]=[a,a^{-i}b]$ has order $2^{n}$. If $c=a^{-i}b$, then
the map $\varphi$ on $H$ defined by
$(a^sc^tx)^{\varphi}=(ac^{2^{n-1}})^sc^tx$ for all $x\in Z(H)$
and integers $s,t$, is a noninner automorphism of $H$ of order 2
leaving $Z(H)$ elementwise fixed. If $a^{2^n}=b^{2^ni}$ and $i$
is even, then we can similarly construct such a $\varphi \in
Aut(H)$.\\
 Hence,  from now on we may assume that $a^{2^ni}=b^{2^n}$ for some
odd integer $i$ and so $c=a^{-i}b$ has order $2^{n+1}$.\\
Now suppose that  $[a,b]\in \left<a^{2^n}\right>$. Then
$Z(H)=\left<a^{2^n}\right>$ and so $|a^{2^n}|\geq 2^n$. Thus
$a^{2^nj}=c^{2^n}$ for some integer $j$. Since $n\geq 2$,
$|a^{2^n}|\geq 2^{n}$ and $|c|=2^{n+1}$, $j$ must be even. This
implies that  $d=a^{-j}c$  has  order $2^n$ and $d^{2^{n-1}}\notin
Z(H)$, as $[a,b]=[a,d]$ is of order $2^n$. Hence the map
$\varphi$ on $G$ defined by
$(a^sd^tx)^{\varphi}=(ad^{2^{n-1}})^sd^tx$ for all $x\in Z(H)$
and integers $s,t$ is the desired automorphism
$\varphi$ of $H$.\\
Thus we may assume that $[a,b]\not\in \left<a^{2^n}\right>$. Since
$Z(H)=\left<a^{2^n},[a,b]\right>$ is cyclic, it follows that
$Z(H)=\left<[a,b]\right>=H'$. On the other hand
$$\frac{H}{Z(H)}=\left<aZ(H)\right>\times \left<bZ(H)\right>$$ and
$|\left<aZ(H)\right>|=|\left<bZ(H)\right>|=2^n$,  which implies
that the element $e=a^{-2^{n-1}i}b^{2^{n-1}}$ does not belong to
$Z(H)$ and $|e|=2$ as $n\geq 2$. Now the map $\varphi$ on $H$
defined by $(a^sb^tx)^{\varphi}=(ae)^s (be)^tx$ for all $x\in
Z(H)$ and integers $s,t$ is the required automorphism $\varphi$.
This completes the proof. \hfill $\Box$\\

\noindent{\bf Acknowledgement.} The author thanks the Center of
Excellence for Mathematics, University of Isfahan.


\end{document}